\documentclass{amsart}
\usepackage{amsmath,amsthm}
\usepackage{color}
\usepackage{epsfig}

\theoremstyle{plain}
\newtheorem{thm}{Theorem}[section]

\newtheorem{claim}[thm]{Claim}

\newtheorem*{thm-eight}{Theorem 8.1}

\theoremstyle{definition}
\newtheorem{defn}[thm]{Definition}
\newtheorem{remark}[thm]{Remark}

\newtheorem*{problem}{Problem 3.6(D)}

\newcommand{\hatR}{\ensuremath{{\widehat{R}}}}

\newcommand{\calL}{\ensuremath{{\mathcal L}}}

\newcommand{\comment}[1]{}

\DeclareMathOperator{\nbhd}{Nbhd}
\newcommand{\N}{\ensuremath{\mathbb{N}}}

\newcommand{\Z}{\ensuremath{\mathbb{Z}}}

\def\smallcoprod{\raise.3ex\hbox{$\,\scriptstyle\coprod\,$}}

\def\eop#1{\hfill\break\rightline{$\square$\ #1}}

\begin{document}
\title{Infinitely many knots admitting the same integer surgery}
\author{J. Luecke}
\address{University of Texas at Austin}
\email{luecke@math.utexas.edu}
\author{J. Osoinach}
\address{University of Dallas}
\email{josoinach@udallas.edu}

\begin{abstract}
The construction of knots via annular twisting has been used to create families of knots yielding the same manifold via Dehn surgery.  Prior examples have all involved Dehn surgery where the surgery slope is an integral multiple of 2.  In this note we prove that for any integer $n$ there exist infinitely many different knots in $S^3$ such that $n$-surgery on those knots yields the same manifold.  In particular, when $|n|=1$ homology spheres arise from these surgeries. 
In addition, when $n \neq 0$ the bridge numbers of the knots constructed tend 
to infinity as the number of twists along the annulus increases.
\end{abstract}



\maketitle

Dehn surgery on knots is a long-standing technique for the construction of 3-manifolds.  While well-known theorems of Lickorish \cite{lickorish} and Wallace \cite{wallace} state that every orientable 3-manifold can be obtained by Dehn surgery on some link in $S^3,$ this representation is far from unique.  In particular, in the Kirby problem list \cite{kirby}, Clark asks the following question:

\begin{problem}
  Is there a homology 3-sphere (or any 3-manifold) which can
be obtained by $n$-surgery on an infinite number of distinct
knots?
\end{problem}

\noindent In \cite{osoinach}, the parenthetical version of this question was answered affirmatively by constructing knots using the method of twisting along an annulus. This method was subsequently developed and refined in \cite{teragaito} to construct knots yielding a small Seifert-fibered manifold; in \cite{teragaito10} to describe properties of a toroidal manifold so constructed; and in \cite{bgl} to analyze the bridge number of the knots constructed.  In these 
examples, however, the surgeries used to construct the manifolds from the knots 
have all had even integer
slopes.  In \cite{osoinach} and \cite{teragaito10} the surgery slope is $0$, and in \cite{kouno},\cite{teragaito} and \cite{bgl} the surgery slopes are 
multiples of $4$. An homology sphere results exactly when $|n|=1$.

This note uses the annular twist construction to create, for each integer $n$, an infinite family of distinct knots in $S^3$ such that $n$-surgery on each knot in the collection yields the same manifold. When $|n|=1$, the resulting 
manifold is an homology sphere thereby 
answering affirmatively Problem 3.6(D) above.   The members of each infinite
family are distinguished by their hyperbolic volume. Alternatively,
at least when $n \neq 0,$ the knots in a family are shown to be different by 
proving that the bridge numbers tend to infinity as the number of twists along the annulus increases. 

The Dehn surgeries on a knot, $k$, 
in the $3$-sphere are parameterized by their surgery slopes. These surgery slopes are described by $p/q \in \mathbb{Q} \cup \infty$, meaning that the slope is a curve that runs $p$ times meridianally and $q$ times longitudinally (using the preferred longitude)
along the boundary of the exterior of $k$. We write $k(p/q)$ for the $p/q$ 
Dehn surgery on $k$. In this notation, an $n$-surgery on $k$ refers to the
integer surgery $k(n/1)$. 

\begin{figure}
\centerline{\epsfysize=3truein\epsfbox{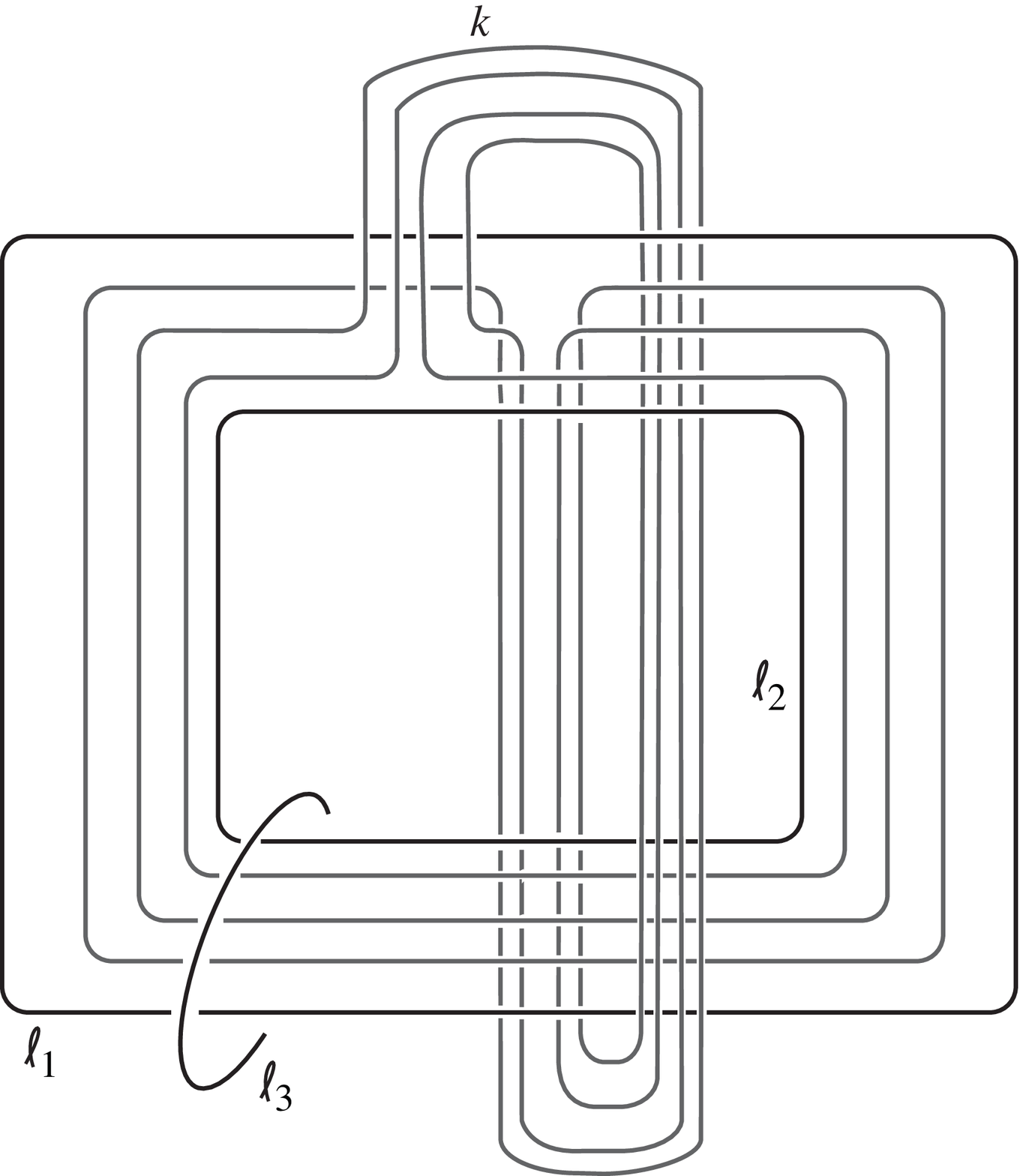}}
\caption{The link $\calL$}
\label{Example}
\end{figure}

\begin{defn}\label{defL}
Let $\calL = k \cup l_1 \cup l_2 \cup l_3$ be the link picture in Figure~\ref{Example}.
Let $L(\alpha, \beta, \delta, \gamma)$ be the corresponding Dehn surgery on $\calL$.
Here the surgery slopes $\alpha,\beta,\delta, \gamma$ will be either in $\mathbb{Q} \cup \infty$,
using the meridian-longitude coordinates on the boundary of a knot in $S^3$ (with a right-handed orientation
on $S^3$), or an asterisk, meaning
that no surgery is done on that component and the component is seen as a knot in the
surgered manifold. 
\end{defn}

The main result of this note is the following theorem.

\begin{thm}\label{mainthm}
For integers $m,n$, $k_n^m=L(*,-1/m,1/m,-1/n)$ is a knot in $S^3$. Furthermore $k_n^{m_1}(n)$
is homeomorphic to $k_n^{m_2}(n)$ for any integers $m_1,m_2$.
\begin{enumerate}
\item  For a fixed $n \neq 0$,
the bridge number of $k_n^m$ tends to infinity as $m$ tends to infinity.
\item For any integer $n$ there
is an $M_n>0$ such that if $m_2 > m_1 > M_n$ then $k_n^{m_1}$ and $k_n^{m_2}$ 
are hyperbolic knots with the hyperbolic volume of $k_n^{m_2}$ larger than
that of $k_n^{m_1}$.
\end{enumerate}
In particular, for each integer $n$ there are infinitely many different knots in the 
family $\{k_n^m\}$.
\end{thm}

\noindent{\it Proof.}

\begin{figure}
\centerline{\epsfysize=3truein\epsfbox{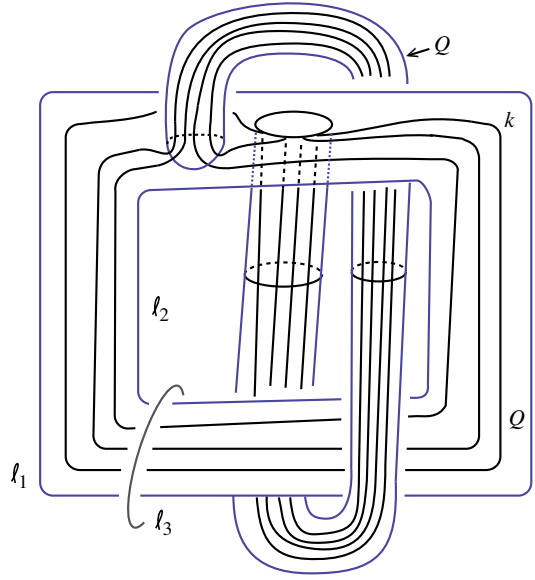}}
\caption{The $2$-punctured Klein bottle $Q$ containing $k$}
\label{Q}
\end{figure}

We first show that for any integer $n$, the $n$-surgery on each $k_n^m$ yields the same
manifold for each $m$.
Figure~\ref{Q} shows that the knot $k$ is a non-separating, orientation-preserving curve on a twice-punctured Klein bottle,
$Q$, cobounded by $l_1$ and $l_2$ and in the complement of $l_3$.

Thus $Q - \nbhd(\calL)$
is a $4$-punctured sphere, $P$, properly embedded in the exterior, $X_{\calL}$, of
$\calL$ in $S^3$. $\partial P$ has one component on each of
$\partial \nbhd(l_1),\partial \nbhd(l_2)$ of slope $0/1$ and
two components on $\partial \nbhd(k)$ of slope $0/1$. Let $\widehat{P}$ be the
annulus in the exterior of
$L(0/1,*,*,*)$ obtained by capping off the two components of $P$ along $\partial \nbhd(k)$.
Dehn twisting this exterior along the properly embedded
annulus $\widehat{P}$ in $L(0/1,*,*,*)$ (see Remark~\ref{twisting}), induces a homeomorphism
of $L(0/1,-1/m,1/m,-1/n)$ to $L(0/1,1/0,1/0,-1/n)$ for each $m,n$.

\begin{remark}\label{twisting}
Let $\hatR$ be an annulus embedded in a $3$-manifold $M$ with $\partial \hatR$ the link
$L_1 \cup L_2$ in $M$. Let $R=\hatR \cap (M - \nbhd(L_1 \cup L_2))$.
Fix an orientation on $M$ and $\hatR$. This induces
an orientation on $L_i$ and its meridian $\mu_i$.
Let $\hatR \times [0,1]$ be a product
neighborhood of $\hatR$ in $M$ so that the corresponding interval orientation on $R \times [0,1]$
corresponds to the meridian orientation of $L_1$. Pick coordinates
$\hatR=e^{2 \pi i \theta} \times [0,1]$, with $\theta \in [0,1]$, so that
$e^{2 \pi i \theta} \times \{0\}, \theta \in [0,1],$ is the oriented $L_1$.
Define the homeomorphism
$f_m\colon  \hatR \times [0,1] \to \hatR \times [0,1]$ by $(e^{2 \pi i \theta},s,t) \to
(e^{2 \pi i(\theta + mt)},s,t)$. Note that $f_n$ restricted to $\hatR \times \{0,1\}$
is the identity. Assume that the knot $K$ in $M$ intersects $\hatR \times
[0,1]$ in $[0,1]$ fibers. Let $K^m$ be the knot in $M$ gotten by applying
$f_m$ to $K
\cap (\hatR \times [0,1])$ (and the identity on $K$ outside this region).
We refer to $K^m$ as {\em $K$ twisted $m$ times
along $\hatR$}, or we say that {\em $K^m$ is obtained from $K$ by twisting
along $\hatR$}. Furthermore, note that $f_m$ induces a homeomorphism
$h_m \colon M-\nbhd(L_1 \cup L_2) \to M-\nbhd(L_1 \cup L_2)$ by applying $f_m$
in $R \times [0,1]$ along with the identity outside this neighborhood.  We refer to this homeomorphism
$h_m$ of $M-\nbhd(L_1 \cup L_2)$  as {\em $m$ Dehn-twisting along the properly embedded annulus $R$}.
\end{remark}

\begin{figure}
\centerline{\epsfysize=1.5truein\epsfbox{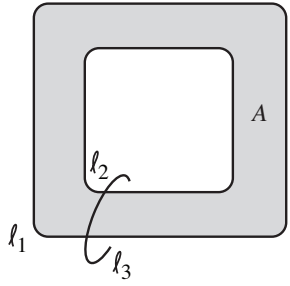}}
\caption{The annulus $A$ bounded by $l_1 \cup l_2$}
\label{A}
\end{figure}

\bigskip

Figure~\ref{A} shows an annulus $A$ cobounded by $l_1$ and $l_2$ in the complement of $l_3$ (which can be taken to
intersect $k$ algebraically zero and geometrically four times and which induces the framing $0/1$ on each of $l_1$ and $l_2$), which becomes an annulus $A_n$ cobounded by $l_1$ and $l_2$ after $(-1/n)$-surgery on $l_3$. Dehn-twisting the
exterior of $l_1 \cup l_2$ in $L(*,*,*, -1/n)$ along $A_n$ (really the restriction of $A_n$ to this exterior, Remark~\ref{twisting}) induces an orientation-preserving homeomorphism
of the manifold $L(1/0,-1/0,1/0,-1/n)=S^3$ to the manifold $L(1/0, -1/m, 1/m, -1/n)$. The inverse of this
homeomorphism identifies $k_n^m$ as a knot in $S^3$ obtained by twisting $k_n^0$ along $A_n$ (see Remark~\ref{twisting}).

The following claim finishes the argument that the $n$-surgeries on $k_n^m$ are the same manifold.

\begin{claim}\label{clm:slope}
For each $m,n$, $L(0/1,-1/m,1/m,-1/n)=k_n^m(n)$.
\end{claim}

\noindent {\it Proof of Claim:} $L(0/1,-1/m,1/m,-1/n)$ is clearly a surgery on $k_n^m$. Our goal is to identify the slope of this surgery, $\alpha(m,n)$, in terms of the coordinates on $k_n^m$ as a knot in $S^3$.
Let $P_n$ be the $4$-punctured sphere $P$ after $-1/n$-surgery on $l_3$. Then $\alpha(0,n)$ is the slope of $P_n$ on $k_n^0$.

Twisting along $A_n$ gives a homeomorphism of the exterior
of $l_1 \cup l_2 \cup k_n^0$ in $S^3$ to the exterior of $l_1 \cup l_2 \cup k_n^m$ and
consequently takes $P_n$ to a $4$-punctured sphere $P_n^m$ in the exterior of $l_1 \cup
l_2 \cup k_n^m$. The slope $\alpha(m,n)$ is the slope of $P_n^m$ on $k_n^m$. We may use
$P_n^m$ to compute the linking number of the slope $\alpha(m,n)$ with $k_n^m$ and consequently the
coordinates of the slope. Orient $k_n^m$ and take the orientation on $P_n^m$ that induces an
orientation on $\partial P_n^m \cap \nbhd(k_n^m)$ that agrees with that on $k_n^m$. Then twice
the linking number of $\alpha(m,n)$ with the oriented $k_n^m$ in $S^3$ is the negative of the linking number between the oriented $k_n^m$ and $l_1 \cup l_2$, given the orientation induced by $P_n^m$ on $l_1 \cup l_2$. By considering $k_n^m$ as twisting $k_n^0$ along $A_n$ away from $l_1 \cup l_2$,
one sees that this latter linking number is $-2n$ (one may verify that in the $1/0$ surgery on $l_3$, this
linking number is zero, then observe how the linking number changes under $-1/n$-surgery). Thus $\alpha(m,n)$ is the slope $n/1$ as desired.

\eop{(Claim~\ref{clm:slope})}

\begin{claim}\label{clm:noannulus}
Let $X_n$ be the exterior of $L(*,*,*,-1/n)$ and $T_1,T_2$ be the components of
$\partial X_n$ coming from $\nbhd(l_1),\nbhd(l_2)$ respectively.
For each integer $n \neq -2$, the interior of $X_n$ is hyperbolic.
For every integer $n$ (including $-2$), there is no essential annulus properly 
embedded in 
$X_n$ with
one boundary component on $T_1$ and the other on $T_2$.
\end{claim}

\begin{figure}
\centerline{\epsfysize=4truein\epsfbox{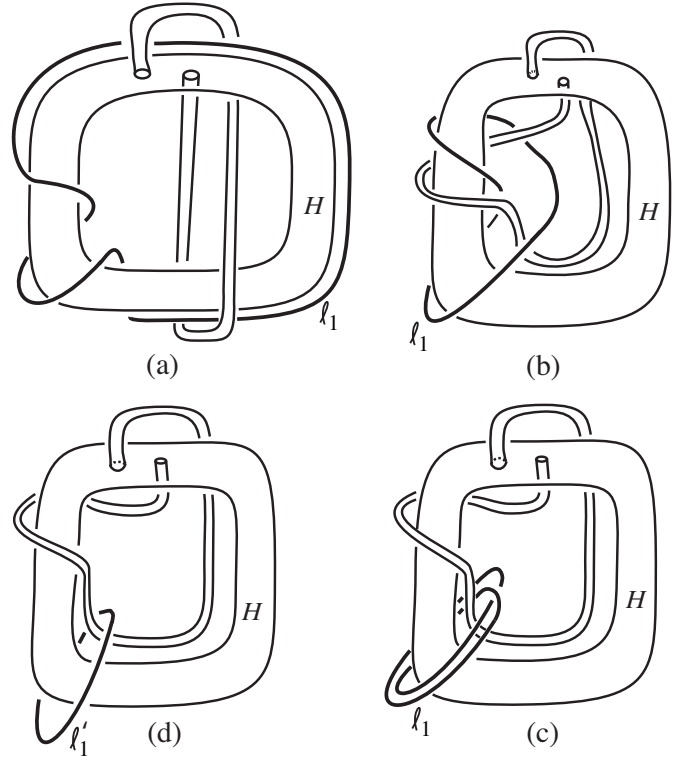}}
\caption{ $l_1$ in  $L(*,*,*,1/2)$ is a $(2,-1)$-cable on $l_1'$}
\label{2surgery}
\end{figure}

\noindent {\it Proof of Claim:}
SnapPy \cite{snappy} shows that $\calL$ is hyperbolic. The program HIKMOT 
\cite{hikmot} certifies this calculation.
The sequence of isotopies Figure~\ref{2surgery}(a)-(c) shows that  $l_1$ in  $L(*,*,*,1/2)$ is a $(2,-1)$-cable on the knot $l_1'$ pictured
in Figure~\ref{2surgery}(d) (the 3-manifold $H$ in Figure~\ref{2surgery} is a
neighborhood of the punctured Klein bottle $Q$ and $l_1$ is pushed
off $H$). Because the linking number of $l_1'$ with $k$ is one, the exterior of $k \cup l_1 \cup l_2$  in
$L(*,*,*,1/2)$ is toroidal. It follows from \cite{gordon} 
and \cite{gordon-wu} that the interior of $X_n$ is
hyperbolic as long as $|n+2|>3$.

For $n \in \{1,0,-1,-3,-4,-5 \}$ SnapPy shows that $X_n$ is hyperbolic 
and HIKMOT certifies this calculation. 
Thus the interior of $X_n$ is hyperbolic, and in particular $X_n$ is annular,
as long as $n \neq -2$.

We must still show that $X_{-2}$ is anannular. As mentioned above,
Figure~\ref{2surgery}(d) shows that $X_{-2}$ is the union, along a
torus $T$, of the exterior of a $(2,-1)$-cable of the core
of a solid torus and
the exterior, $X_{-2}'$, of $l_1' \cup l_2 \cup k$ after $1/2$-surgery on $l_3$.
SnapPy shows $X_{-2}'$ is hyperbolic and HIKMOT certifies this. Now assume there
were an essential annulus in $X_{-2}$ between $T_1$ and $T_2$, and consider its
intersection with the incompressible torus $T$. We may surger away any closed curves of intersection
which are trivial on $T$. Then an outermost component of intersection
with $X_2'$ will give rise to an essential annulus or disk
properly embedded in $X_2'$, contradicting the hyperbolicity of $X_2'$.
\eop{(Claim~\ref{clm:noannulus})}

We first verify $(1)$ of Theorem~\ref{mainthm}.
As before, let $A_n$ be the annulus from Figure~\ref{A} cobounded by $l_1$ and $l_2$  and after $(-1/n)$-surgery on $l_3$. The
knot $k_n^m$ is obtained by twisting $k$  along $A_n$ ($m$ times) in the copy of $S^3$ obtained by $(-1/n)$-surgery
on $l_3$.  As the linking number of $l_1$ and $l_2$ in this copy of $S^3$ is $n$, $l_1 \cup l_2$ is not
the trivial link. Then Claim~\ref{clm:noannulus} along with Corollary 1.4 of \cite{bgl} shows that for $n \neq 0$ the (genus $0$) bridge number of the knots $k_n^m$ in $S^3$ goes to infinity  as $m$ goes to infinity
(as the linking number of $l_1$ and $l_2$ is non-zero, Lemma 2.4 of \cite{bgl} shows there is a catching surface for the
pair $(k, A_n)$).
Note that since $A_0$ lies on a Heegaard sphere for $S^3$, 
the bridge numbers of $\{k_0^m\}$ will be bounded.

We now verify $(2)$ of Theorem~\ref{mainthm}. By Claim~\ref{clm:noannulus},
the interior of $X_n$ is hyperbolic whenever $n \neq -2$. 
Thurston's Dehn Surgery Theorem and Theorem 1A of \cite{neumann},
shows that there is an $M_n>0$
such that for $m >M_n$, $k_n^m$ is hyperbolic and its volume increases monotonically with
$m$. When $n = -2$, recall from the proof of Claim~\ref{clm:noannulus} that
Figure~\ref{2surgery}(d) shows that $X_{-2}$ is the union, along a
torus $T$, of the exterior of a $(2,-1)$-cable of the core
of a solid torus and
the exterior, $X_{-2}'$, of $l_1' \cup l_2 \cup k$ after $1/2$-surgery on $l_3$.
That is, identify $L(*,*,*,1/2)$ as a link in $S^3$ by putting two full 
left-handed twists along the linking circle $l_3$. Then $L(*, -1/m,1/m,1/2)$
corresponds to $(-1-2m)/m$ surgery on $l_1$ and $(1-2m)/m$ surgery on $l_2$.
The Seifert fiber on $l_1$ as a $(2,-1)$-cabling on $l_1'$ is $-2/1$. As the 
surgery slope intersects this Seifert fiber slope once, this surgery on
$l_2$ corresponds to doing a $(-1-2m)/4m$ surgery  on $l_1'$ 
(see Corollary 7.3 of \cite{gordon2}). As noted above,
HIKMOT verifies $k \cup l_1' \cup l_2$ to be hyperbolic.
Thus an application of Theorem 1A of \cite{neumann} to the exterior $X_{-2}$
of this link, shows there is an $M_{-2}$ such that for $m > M_{-2}$, $k_{-2}^m$
is hyperbolic and its volume increases monotonically with $m$.

Since hyperbolic volume and bridge number are knot invariants, either $(1)$
 (when $n \neq 0$)
or $(2)$ shows that for an integer $n$ the family $\{k_n^m\}$ is infinite.
\eop{(Theorem~\ref{mainthm})}

\begin{remark} SnapPy shows the homology spheres that arise 
in the above construction ($|n|=1$) to be hyperbolic manifolds with 
$volume(k_{-1}^0(-1))=3.400436870$ and $volume(k_{1}^0(1))=5.7167678901$. 
SnapPy shows the manifold corresponding
to $n=-2$ to be hyperbolic with $volume(k_{-2}^0(-2))=3.110698158$.
These calculations are not verified by HIKMOT.
\end{remark}

When $K$ is a knot in the $3$-sphere, let $W(K,n)$ be the $4$-manifold obtained
by attaching a 2-handle to the $4$-ball along $S^3$ with framing $n$.
The 3-manifold $K(n)$
is the boundary of $W(K,n)$. In \cite{ajot}, the following analog of
Problem 1 is considered.

\bigskip

\noindent {\bf Problem 2} {\it Let $n$ be an integer. Find infinitely many
mutually distinct knots $K_1, K_2, \dots$ such that $W(K_i,n)$ is
diffeomorphic to $W(K_j,n)$ for each $i, j \in \N$.}

\bigskip

In \cite{ajot}, such an infinite family of knots is demonstrated when $n \in
\{-4,0,4\}$. These families of knots are constructed by twisting along
an annulus similar to the construction above for the family
$\{k_n^m, m \in \Z \}$.  Indeed
Theorem~\ref{mainthm} shows that the boundaries of $W(k_n^{i},n)$ and
$W(k_n^{j},n)$ are diffeomorphic. So it is natural to ask

\bigskip

\noindent {\bf Problem 3} {\it Let $n$ be an integer. Are $W(k_n^{i},n)$ and
$W(k_n^{j},n)$ diffeomorphic?}

\bigskip

\noindent{\bf Acknowledgements.}
The authors would like to thank Kyle Larson for very helpful conversations,
and Neil Hoffman for his help with HIKMOT.

\providecommand{\bysame}{\leavevmode\hbox to3em{\hrulefill}\thinspace}
\providecommand{\MR}{\relax\ifhmode\unskip\space\fi MR }
\providecommand{\MRhref}[2]{%
  \href{http://www.ams.org/mathscinet-getitem?mr=#1}{#2}
}
\providecommand{\href}[2]{#2}

\end{document}